\documentclass[12pt,leqno,twoside]{amsart}
\usepackage{amssymb,amsfonts,amsthm,amsmath,amscd}
\usepackage{latexsym}

\newtheorem{theorem}{Theorem}[section]
\newtheorem{corollary}[theorem]{Corollary}
\newtheorem{lemma}[theorem]{Lemma}
\newtheorem{proposition}[theorem]{Proposition}
\theoremstyle{definition}
\newtheorem{definition}[theorem]{Definition}
\theoremstyle{remark}

\newtheorem{note}[theorem]{\sc Note}
\renewcommand{\Box}{\square}    %\diamond

\newcommand{\cod}{{\rm{cod}}}

\newcommand{\reg}{{\rm{reg}}}

\newcommand{\G}{{\rm{G}}}

\newcommand{\Eu}{{\mathrm{Eu}}}
\newcommand{\sing}{{\rm{sing}}}

\newcommand{\rank}{{\rm{rank\hspace{2pt}}}}
\newcommand{\Crit}{{\rm{Crit\hspace{2pt}}}}

\newcommand{\LMD}{{\rm{LMD}}}
\newcommand{\NMD}{{\rm{NMD}}}

\newcommand{\M}{{\rm{M}}}
\newcommand{\Ma}{{\rm{Ma}}}

\newcommand{\cl}{{\rm{closure}}}

\newcommand{\ity}{{\infty}}

\newcommand{\supp}{{\rm{supp}}}

\newcommand{\fin}{\hspace*{\fill}$\Box$\vspace*{2mm}}

\newcommand{\cO}{{\mathcal O}}

\newcommand{\cS}{{\mathcal S}}

\newcommand{\bC}{{\mathbb C}}
\newcommand{\bP}{{\mathbb P}}

\newcommand{\bZ}{{\mathbb Z}}

\begin{document}

\title[Index formula for MacPherson cycles of affine algebraic varieties]
{Index formula for MacPherson cycles of affine algebraic varieties}

\author{\sc J\"org Sch\"urmann}

\address{Westf\" alische Wilhelms-Universit\" at \endgraf
Mathematisches Institut \endgraf
  Einsteinstr. 62 \endgraf
  48149 M\" unster \endgraf
 Germany}

\email{jschuerm@math.uni-muenster.de}

\author{\sc Mihai Tib\u ar}

\address{Math\' ematiques \endgraf
Universit\'e de Lille 1 \endgraf   59655 Villeneuve d'Ascq \endgraf France}

\email{tibar@math.univ-lille1.fr}

%\date{October 24, 2006}

\subjclass[2000]{Primary 14C25; Secondary 14C17, 14R25, 32S60, 14D06, 32S20}

\keywords{characteristic classes, constructible function, affine polar varieties, 
Euler obstruction, index theorem, characteristic cycles, stratified Morse theory}

\begin{abstract}
We give explicit MacPherson cycles
for the Chern-MacPherson class of a closed affine algebraic variety $X$ and for any 
constructible function $\alpha$ with respect to a complex algebraic Whitney stratification of $X$. 

We define generalized degrees of the global polar varieties and of the MacPherson cycles 
and we prove a global index formula for the Euler characteristic
of $\alpha$.
  Whenever $\alpha$ is the Euler obstruction of $X$, this index formula specializes
 to the Seade-Tib\u ar-Verjovsky global counterpart of the L\^e-Teissier formula for the local Euler obstruction. 
\end{abstract}
%
%%% ---------------------------------------------------------------------
\maketitle

%%% ----------------------------------------------------------------------
%\tableofcontents
%%% ----------------------------------------------------------------------

\setcounter{section}{0}
\section{Introduction}\label{intro}
%%%%%%%%%%%%%%%%%%%%%%%%%%%%%%

 MacPherson \cite{MP} has defined analogues of Chern classes for singular algebraic (or analytic) varieties $X$, based on Chern-Mather classes
 and his famous  \emph{local Euler obstruction} $\Eu_X$, a constructible function on $X$ defined by obstruction theory on the Nash blow-up, which measures in some sense the singularities of $X$.
  It is a natural challenge to find representing cycles for these Chern-MacPherson homology classes, such that they reflect the geometry of $X$ up to a certain extent. In the local analytic case, L\^e and Teissier \cite{LT} explained how the generic local polar varieties enter in the description of Chern-MacPherson classes,
  in particular they expressed the local Euler obstruction at a given point as the alternating sum of suitable polar multiplicites.  It turned out that these polar varieties are the ingredients of the desired representing cycles. Later on this was worked out by Massey \cite{Ma1} in greater generality,
  based on a reformulation of MacPherson's theory in terms of characteristic cycles of constructible sheaves or holonomic $D$-modules, see also \cite{BDK, Du2, Gi, Sa, Sch-lect}.

We develop here the global affine algebraic counterpart. Abstractly, it is known that the Chern-MacPherson classes of a complex algebraic proper subset $X \subset \bC^N$ can be represented by algebraic cycles since the MacPherson transformation may be defined by using Chow groups, see e.g. \cite{Ken}. Our aim is to produce an explicit \emph{global geometric MacPherson cycle}. This construction needs techniques adapted to the affine global setting. We therefore use global general coordinates in order to introduce the key new tool: \emph{affine polar varieties}. We give in the following a brief account of our results.

%%%%%%%%%%%%%%%%%%%%%%%%%%%%%%
Let us assume for the moment that $X$ is of pure dimension $n<N$. The \emph{$k$-th global polar variety of $X$} ($0\leq k \leq n$) is the following algebraic set:
\[  P_k(X) = \overline{\Crit (x_1, \ldots , x_{k+1})_{|X_\reg} }, \]
with $\Crit (x_1, \ldots , x_{k+1})_{|X_\reg}$ the usual critical locus of points $x\in X_{\reg}$
where the differentials of these functions restricted to $X_{\reg}$ are linearly dependent.
 For \emph{general} coordinates $x_i$, the polar variety $P_k(X)$ has pure dimension $k$ or it is empty, for all $0\le k< n$. We have $P_n(X) := X$ 
 and  we set $P_k(X) := \emptyset$ for $k>n$.
It turns out that $P_k(X)$ represents
the \emph{$k$-th dual Chern-Mather class} of $X$. Note that we index our polar varieties by their dimension, and not by their codimension as often done.

We fix an algebraic Whitney stratification $\cS$ with connected strata.
In this context $X$ need not be pure dimensional and we only assume $n=\dim X<N$.
 Let $\alpha : X\to \bZ$ be an $\cS$-constructible function, meaning that the restriction $\alpha_{|S}$ is locally constant for all strata $S\in \cS$.
By stratified Morse theory \cite{GM, Sch-book},  one attaches to every $S$ the normal Morse index $\eta(S, \alpha)$,
which is the Euler characteristic of the normal Morse data $\NMD(S)$ of $S$ weighted by $\alpha$. 

We define the \emph{$k$-th MacPherson cycle of $\alpha$} ($0\leq k \leq n$) by:
\begin{equation}\label{eq:macph}
\Lambda_k(\alpha) = \sum_{S\in \cS} (-1)^{\dim S} \eta(S,\alpha) P_k(\bar S),
\end{equation}
where $P_k(\bar S)$ is the $k$-th global polar variety of the algebraic closure 
$\bar S\subset \bC^N$ of the stratum $S$. Here we assume that
the coordinates $\{x_i\}$ are chosen \emph{general} with respect to all
 strata $S$ of our algebraic Whitney stratification.
The precise meaning of ``general'' will be explained in section \ref{s:MPcycles}.
We then use a reinterpretation of the MacPherson transformation (cf. \S \ref{ss:macph}) in order to prove the following statement, to which Theorem \ref{t:MPcycle2} represents a slight extension:
\begin{theorem}\label{t:MPcycle}
 
For any $\cS$-constructible function $\alpha$, the cycle class $[\Lambda_k(\alpha)]$ of $\Lambda_k(\alpha)$ represents the $k$-th dual 
Chern-MacPherson class $\check c_k^M(\alpha)$ in the Borel-Moore homology group $H^{BM}_{2k}(X)$
or in the Chow group $CH_k(X)$. 
\end{theorem}

Our main result is a general \emph{index formula} for constructible functions. We explain how it is related to other formulas in the literature.
To get the degrees $\gamma_k(X)$ of our polar varieties $P_k(X)$
and the generalized degrees $\gamma_k(\alpha)$ of our MacPherson cycles $\Lambda_k(\alpha)$,
we have to use additional generality conditions which take into consideration the hyperplane at infinity of the projective closure $X\subset \bar X \subset \bP^N$.
  This generality of the system of coordinates $\{ x_i\}_i$ implies 
  in case of a pure dimensional $X$ (see Proposition \ref{p:degree} and compare with 
\cite{Ha}) that the projection $(x_1, \ldots , x_{k}) : P_k(X) \to \bC^k$
is proper and finite for $1\leq k \leq n$, so that its degree $\gamma_k(X)$ is well defined.
We may also define $\gamma_0(X) := \# P_0(X)$ and for $k>n$ it is natural to set $\gamma_k(X) = 0$. 
  Once the affine coordinates of $\bC^N$ are fixed, these degrees are independent on linear change of coordinates, provided general. Nevertheless they depend on the embedding of $X$ in $\bC^N$.

 The degrees occur in the following global index formula:
\begin{equation}\label{eq:indexeu}
\Eu(X) = \sum_{k=0}^n (-1)^{n-k}\gamma_k(X) ,
\end{equation} 
where $\Eu(X)$ is the \emph{global Euler obstruction} introduced in \cite{STV}. This formula, proved in \cite{STV} by using stratified vector fields techniques, is the global counterpart of the local formula by L\^e-Teissier \cite{LT}. A
more geometric interpretation for $\gamma_k(X)$ is the number of complex Morse points of the affine pencil defined by the general coordinate $x_{k+1}$ on the \emph{regular part} of the general slice $X\cap \{ x_1=t_1, \ldots, x_k=t_k\}$; by definition, this slice is $X$ for $k=0$,
and $\emptyset$ for $k\ge n+1$.\\
% In terms of global Euler obstructions, this key ingredient can be 
% also written as follows, cf. \cite{STV}:
% 
% \[  \gamma_k(X) = (-1)^{n-k} [ \Eu(X\cap \{ x_1=t_1, \ldots, x_k=t_k\}) - \Eu(X\cap \{ %x_1=t_1, \ldots, x_{k+1}=t_{k+1}\})].\]

More generally, we define for some
$\cS$-constructible function $\alpha$
the \emph{(generalised) degree} $\gamma_k(\alpha)$ of $\Lambda_k(\alpha)$ as
\[ \gamma_k(\alpha) := \sum_{S\in \cS} (-1)^{\dim S} \eta(S,\alpha)\gamma_k(\bar S)\:.\]
The definition of $\gamma_k(\alpha)$ for $k>0$ only works in the affine algebraic context and is a \emph{global} counterpart of
the local intersection number 
\[\gamma_k(\alpha)(p):= \sharp_p \;([\Lambda_k(\alpha)] \cap [\{x_1=p_1, \ldots, x_k=p_k\}] )\]
at a point $p=(p_1,\dots,p_N)\in  \bC^N$ as used by Massey \cite[\S 4]{Ma1} in his definition of 
\emph{characteristic polar multiplicities}.

We recall that the Euler characteristic weighted by a constructible function $\alpha$ is defined as follows:
\[ \chi(X,\alpha) := \sum_{S\in \cS} \alpha(S)\cdot \chi(H^*_c(S)),\]
where the compact support cohomology $H^*_c(S)$ is finite dimensional since the strata are locally closed algebraic sets.
We then prove
the following index theorem: 
\begin{theorem}\label{t:mainindex}
For general coordinates $\{x_i\}_{i=1}^N$ and for any $\cS$-constructible function $\alpha$ and $0\leq k \leq n$, we have
\begin{equation}\label{eq:mainindex}
(-1)^k\gamma_k(\alpha) = \chi(X\cap \{x_1=t_1, \ldots, x_k=t_k\},\alpha) - \chi(X\cap \{x_1=t_1, \ldots, x_{k+1}=t_{k+1}\},\alpha).
\end{equation}
In particular,
\begin{equation}\label{eq:indexchi}
\chi(X,\alpha) = \sum_{k=0}^n (-1)^k \gamma_k(\alpha).
\end{equation}
\end{theorem}
We need to carefully define what exactly means ``general coordinates'' in our affine
context, and this is done in \S \ref{ss:deg}.

It turns out that $\chi(X\cap \{x_1=t_1, \ldots, x_k=t_k\},\alpha)$ is independent of the choice of the coordinates $\{x_i\}$
and values $\{t_i\}$, provided they are chosen general enough. So (\ref{eq:mainindex}) implies the same property also for $\gamma_k(\alpha)$.
For $\alpha=1_X$, formula (\ref{eq:indexchi}) calculates the global Euler characteristic. Such a formula for $\chi(X)$ has been first proved for hypersurfaces with isolated singularities in \cite{Ti-equi} and then for general $X$ in \cite{Ti-Dual}.

For pure dimensional $X$ we may also take $\alpha = \Eu_X$ and then formula (\ref{eq:indexchi})
 specializes to the index formula (\ref{eq:indexeu}).

%In this way many of the results of \cite{Ma1,Ma2,Ma3} can also be developed in the global %affine algebraic context
%(e.g. ``Morse inequalities, graded enriched cycles or monodromy''). But for simplicity we %restrict ourself in this paper to the context of
%constructible functions.

%%%%%%%%%%%%%%%%%%%%%%
%%%%%%%%%%%%%%%%%%%%%%%%%%%%%%%%%
\section{MacPherson's Chern class transformation and characteristic cycles}\label{macph+char}
\subsection{Reformulation of the MacPherson transformation}\label{ss:macph}
Let us recall the main ingredients in MacPherson's definition of his
dual Chern classes of a constructible function and the well known by now
relation to the \emph{theory of characteristic cycles}
(cf. \cite{BDK, Du2, Gi, Ken, Sa, Sch-lect}).
Here we work in the embedded complex analytic or algebraic
context, with $X$ a closed subspace of positive codimension
in the complex manifold $M$. Then the main characters
of this story can be best visualized in the commutative diagram

\begin{equation} \label{eq:conormal}
\begin{CD}
F(X) @< \check{E}u < \sim < Z(X) @> \check{c}^{\Ma}_{*} >> H_{*}(X) \\
@|  @V cn V \wr V @| \\
F(X) @ > CC > \sim > L(X,M) @> c^{*}(T^{*}M|X)\cap s_{*} >> H_{*}(X) .
\end{CD}
\end{equation}
\noindent
Here $F(X)$ and $Z(X)$ are the groups of constructible functions and cycles
in the corresponding complex analytic or algebraic context. Similarly
$H_{*}(X)$ is either the Borel-Moore homology group in even degrees $H^{BM}_{2*}(X,\bZ)$
or the Chow group $CH_{*}(X)$. \\

The transformation $\check{E}u$ associates to an irreducible subset $Z$ of $X$
the constructible function $\check{E}u_{Z}:=(-1)^{\dim(Z)}\cdot \Eu_{Z}$,
and is linearly extended to cycles. Then $\check{E}u$ is an isomorphism of groups,
since $\Eu_{Z}|Z_{\reg}$ is constant of value $1$. The transformation
$\check{c}^{\Ma}_{*}$ is similarly defined by associating to an irreducible $Z$
the total dual Chern-Mather class $\check{c}^{\Ma}_{*}(Z)$ of $Z$ viewed
in the homology $H_{*}(X)$ of $X$. 
One has the following
description of the dual Chern-Mather class of $Z$ 
in terms of the Segre class of the 
\emph{conormal space} $T_{Z}^{*}M:=\cl(T_{Z_{\reg}}^{*}M)$ of $Z$ in $M$, which is 
 a conic Lagrangian cycle in $T^{*}M|X$:
\begin{equation}\label{cor:dualMather}
\check{c}^{\Ma}_{*}(Z)=c^{*}(T^{*}M|Z)\cap s_{*}(T_{Z}^{*}M) 
\end{equation}
 after e.g. \cite[Lemme (1.2.1)]{Sa} or \cite[Lemma 1]{Ken}.
 The Segre class is defined by 
\begin{equation}\label{eq:segre}
s_{*}(T_{Z}^{*}M):=\hat{\pi}'_{*}(c^{*}(\cO(-1))^{-1}\cap [\bP(T_{Z}^{*}M)])
= \sum_{i\geq 0}\: \pi'_{*}(c^{1}(\cO(1))^{i}\cap [\bP(T_{Z}^{*}M)]),
\end{equation}
cf. \cite[Example 4.1.2]{Fu}.
Here $\cO(-1)$ denotes the tautological line subbundle on the projectivisation
$\hat{\pi}':  \bP(T^{*}M|_Z)\to Z$
with $\cO(1)$ as its dual.
 
  By definition,
$L(X,M)$ is the group of all cycles generated by the conormal spaces $T_{Z}^{*}M$.
The vertical map $cn$ in diagram (\ref{eq:conormal}) is  the correspondence $Z\mapsto T_{Z}^{*}M$. Then (\ref{cor:dualMather}) obviously implies the commutativity of the right square in (\ref{eq:conormal}).

The \emph{dual MacPherson Chern class transformation} is defined by
\begin{equation}
\check{c}^{\M}_{*}:=\check{c}^{\Ma}_{*}\circ \check{E}u^{-1}:\: 
F(X)\to H_{*}(X).
\end{equation}
This agrees up to a sign with MacPherson's original definition 
of his Chern class transformation $c^{\M}_{*}$, namely
\begin{equation}
\check{c}^{\M}_{i}(\alpha) = (-1)^{i}\cdot c^{\M}_{i}(\alpha) \in H_{i}(X) .
\end{equation}
As for the Chern-Mather classes, we have the relation
\[
 \check{c}^{\Ma}_i(Z) = (-1)^{dim(Z)+i} \cdot c^{\Ma}_i(Z) \in H_i(Z)
\]
for $Z\subset X$ irreducible.

\subsection{Stratified Morse theory and characteristic cycles}
Let $X$ be a closed complex analytic (resp. algebraic)
subset of the complex (algebraic) manifold $M$, which is endowed with a complex
(algebraic) Whitney stratification $\cS$ with connected strata $S$. One of the main
ingredients in the \emph{stratified Morse theory} of Goresky-MacPherson \cite{GM} in this
complex context is the complex link $l_{S}$ of $X$ attached to each stratum $S$,
which is defined as follows (cf. \cite[p.161, def.2.2]{GM}):

\begin{definition}\label{def:NMD}
Let $x\in S$ and consider a holomorphic function germ $g: (M,x)\to (\bC,0)$ with $dg_{x}\in T_{S}^{*}M$  \emph{normally nondegenerate} (or in shorter terms ``good'') in the sense of
\cite[p.128, def.12.1, and p.160]{GM}, i.e.,  $dg_{x}$ does not vanish on any generalized
complex tangent space $\tau= \lim_{x_{n}\to x}\; T_{x_{n}}S'$ for a sequence
 of points $x_{n}$ in another stratum $S' \neq S$ converging to $x$.
In terms of conormal spaces, this just means $dg_{x}\notin \overline{T_{S'}^{*}M}$
for all strata $S' \neq S$. Take a \emph{normal slice} $N$ to $S$ in $x$, i.e. a germ of
a closed complex submanifold $N$ of $M$ which is transversal to $S$, with
$N\cap S=\{x\}$. Then $g_{|X\cap N}$ has in $x$ an isolated stratified critical 
point with respect to the induced Whitney stratification of  $X\cap N$.
So it defines a local Milnor fibration \cite[p.165, prop.2.4(a)]{GM} with Milnor fibre
$$l_{S}:=X\cap N \cap B_{\delta}(x) \cap \{g=w\} \quad \text{for} \quad
0< |w| \ll \delta \ll 1 \:, $$
which is called a \emph{complex link} of $S$.
Here $B_{\delta}(x)$ is a closed ball of radius $\delta$ in some local coordinates.
Similarly the pair of spaces
$$\NMD(S):=(X\cap N \cap B_{\delta}(x),
X\cap N \cap B_{\delta}(x) \cap \{g=w\}) $$ 
for $0< |w| \ll \delta \ll 1$,
is called the \emph{normal Morse datum} of $S$.
\end{definition}
 In case $g_{|X}$ has an isolated stratified singularity at $x$,  one does not need to take any normal slice $N$ in the definition of the normal Morse data. Then this is called \emph{local Morse datum} and denoted by $\LMD(g,x)$ (cf. also \cite{Le-Oslo, Le}).
  
One of the main results of 
Goresky-MacPherson\cite[p.163, thm.2.3]{GM} tells
us that the stratified homeomorphy type of  $l_{S}$ and $\NMD(S)$
are independent of all choices, and thus these are invariants of the stratum $S$ in $X$.

\begin{definition}\label{def:NMindex}
Let the function $\alpha: X\to \bZ$ be constructible with respect to our stratification
of $X$. Then its \emph{normal Morse index} $\eta(S,\alpha)$ along $S$ is defined as:
\begin{equation}\label{eg:NMD1}
\eta(S,\alpha):=\chi(\NMD(S),\alpha):=\chi(X\cap N \cap B_{\delta}(x),\alpha)
- \chi(l_{S},\alpha) \:.
\end{equation}
\end{definition}

By the discussion before, this is a well defined invariant of the stratum $S$.
That the invariant (\ref{eg:NMD1}) is well defined also follows, by independent arguments,
from \cite[sec.5.0.2]{Sch-book}. In \emph{loc.cit.} are developed the corresponding results for constructible sheaves.

 Note that, by the \emph{local conic structure} of
$X\cap N \cap B_{\delta}(x)$ (see e.g. \cite{BV} and \cite[p.165, prop.2.4(b)]{GM}),
we also have $\chi(X\cap N \cap B_{\delta}(x),\alpha)=
\alpha(x)$ so that
\begin{equation}\label{eg:NMD2}
\eta(S,\alpha)=\chi(\NMD(S),\alpha)= \alpha(x)- \chi(l_{S},\alpha) .
\end{equation}
In particular,
\begin{equation}\label{eg:NMD3}
\eta(S,\alpha)=\chi(\NMD(S),\alpha)= \alpha(x)
\quad \text{for} \quad \dim(S)=\dim(X)\:,
\end{equation}
since $l_{S}=\emptyset$ in this case.\\

\begin{definition} \label{def:CC}
Let $F(\cS)$ be the group of constructible functions for the stratification
$\cS$, and let $L(\cS)$ be the group of \emph{conic Lagrangian cycles} generated
by the conormal spaces $\overline{T^{*}_{S}M}=:T^{*}_{\bar{S}}M$ to the (closures of)
strata $S$ of our stratification. Then 
$$CC:\; F(\cS)\to L(\cS);\ \ \  CC(\alpha)= \sum_{S\in \cS} (-1)^{\dim S} \eta(S,\alpha)
\cdot T^{*}_{\bar{S}}M \:.$$
\end{definition}
Note that in the analytic context we consider locally finite sums.
We get from (\ref{eg:NMD3}) that $CC$ is injective and, by induction on $\dim(X)$,
we get in the same way that $CC$ is also surjective, so that we have the isomorphism:
\begin{equation} \label{eq:refined}
CC:\; F(\cS) \stackrel{\sim}{\to} L(\cS) \:.
\end{equation}
\begin{lemma}\label{lem:CCindependent}
Considered as a map from $F(\cS)$ to $L(X,M)$,  the map $CC$ does not depend on the choice of the stratification $\cS$. \fin
\end{lemma}
For a proof of this simple fact, we refer to  \cite[the discusssion before diagram (5.63), p.327]{Sch-book}.
We get in this way the isomorphism
 $CC:\; F(X) \stackrel{\sim}{\to} L(X,M)$ from the diagram
(\ref{eq:conormal}). Nevertheless, the refined isomorphism (\ref{eq:refined})
contains more information since it involves the fixed complex
Whitney stratification $\cS$. The key role of the (dual) Euler obstruction comes from the fact that the (dual) Euler obstruction $\check{E}u_{\bar{S}}$ is \emph{constructible}
with respect to  $\cS$, and its characteristic cycle is: 
\begin{equation}\label{eq:EuCC}
CC(\check{E}u_{\bar{S}})= T^{*}_{\bar{S}}M \:.
\end{equation}
A quick detailed proof of this basic fact is given in \cite[(5.35), p.293, p.323-324]{Sch-book}.
This result implies the commutativity of the left square in diagram (\ref{eq:conormal}).
We end this section by stating a result which will be used in the proof of the index 
theorem, i.e., Theorem \ref{t:mainindex}. Its proof is a straightforward application of Goresky and MacPherson's stratified Morse theory \cite{GM}. 
Another independent proof follows from \cite[Thm.5.4.1, p.367-368]{Sch-book},
where a much deeper counterpart for constructible sheaves is given. 
\begin{lemma}\label{thm:indmain}
Consider a holomorphic function germ $h: (M,x)\to (\bC,0)$ having at
$x\in S$ a stratified Morse critical point of $h$
with respect to $\cS$. Then, for any $\alpha\in F(\cS)$, 
one has the following relation between the normal and the local 
 Morse data:
\begin{equation} \label{eq:LMD}
\eta(S,\alpha) =(-1)^{\dim S}\cdot \chi(\LMD(h,x),\alpha) .
\end{equation}
\fin
\end{lemma}

%%%%%%%%%%%%%%%%%%%%%%%%%%%%%%%%%%%%%%%%%%%%%%%%
%%%%%%%%%%%%%%%%%%%%%%%%%%
 %%%%%%%%%%%%%%%%%%%%%%%%%%%%%%%%%%%%%%%%%%%%%%%%
%%%%%%%%%%%%%%%%%%

\subsection{Affine polar varieties and MacPherson cycles}\label{polar}\label{s:MPcycles}
 In order to express
 dual Chern Mather classes in terms of polar varieties,
one starts from a ``general transversality''
result of Kleiman \cite{Kl}. This was used in the local analytic context by Teissier  \cite{Te-pol}
for establishing the existence and the main properties of his ``generic local polar varieties''.
The corresponding global result for projective varieties is due to Piene \cite{Pi,Pi2}.
 One may quickly derive from \cite[Theorem 3]{Pi2} our affine statement Proposition \ref{cor:dualMather2} below using the pull-back by the inclusion map $\bC^N \hookrightarrow \bP^N$. We shall give here a proof in the language of conormal spaces, which is naturally related to characteristic cycles. The main reason of our different proof is that this extends to constructible functions and allows us to easily derive Theorem \ref{t:MPcycle2}.

 All the results of \cite{LT,Pi,Pi2,Te-pol} were first proved using the Nash blow-up. For the translation between these two viewpoints of polar sets, we refer to  \cite{HM,Me,Te-pol}.

\begin{definition}\label{d:polar}
Let $\omega := (\omega_1, \ldots , \omega_{k+1})\in \check{\bC}^N$.
The algebraic set 
\[ \begin{array}{c}
P_{k}(X,\omega):=   
 \overline{ \Crit (\omega_1, \ldots , \omega_{k+1})_{X_{\reg}} }
\end{array} \]
is called the \emph{polar locus of $\omega$}. 
\end{definition}
We shall use this definition for linearly independent forms $\omega_1, \ldots , \omega_N$. This agrees with the definition of $P_k(X)$ in the Introduction if we trivialize $T^*\bC^N$ by
\[(\omega_1, \ldots , \omega_N):=(dx_1,  \ldots , dx_N )\]
for general global affine coordinates  $x_i$ on $\bC^N$.\\ 
 
 \begin{proposition}\label{cor:dualMather2}\cite[Theorem 3]{Pi2}
Let $(\omega_1, \ldots , \omega_{n+1})$ be linearly independent vectors in $\check{\bC}^N$
such that the subspace $\check{H}^{N-1-k}$ generated by $(\omega_1, \ldots , \omega_{k+1})$ 
belongs to  the generic set $\Omega^{N-1-k}(X)$ as in Proposition \ref{p:polar} for all
$0\leq k \leq n$. Then, for any such $k$,  $P_{k}\bigl(X,(\omega_1, \ldots , \omega_{k+1})\bigr)$
is pure $k$-dimensional (or empty) with 
\begin{equation} \label{eq:M=P}
\check{c}^{\Ma}_{k}(X)= [P_{k}\bigl(X,(\omega_1, \ldots , \omega_{k+1})\bigr)] \in CH_{k}(X) .
\end{equation} 
\end{proposition}

\begin{proof}

Since  Kleiman's result \cite{Kl} is proved in the algebraic context, its arguments
also apply to our {\em global affine algebraic} setting, as we explain here in a somewhat larger generality.

Assume that $M$ is a complex algebraic manifold with a trivial (co)tangent bundle
of pure dimension $N$, with
\[p: \bP(T^{*}M)= M \times \bP(\check{\bC}^N) \to \bP(\check{\bC}^N)\]
the projection on the last factor. Let $X\subset M$ be a closed complex  algebraic subset of pure dimension $n<N$,
with $T^*_X M=\cl(T_{X_{\reg}}M) \subset T^*M$ its conormal space. Since this is a conic subset, one can look at the
projectivisation $\bP(T^*_X M)\subset \bP(T^{*}M)$. Consider in addition a closed algebraic subset $Z\subset X$,
containing $X_{\sing}$ such that $U:= X\backslash Z \subset X_{\reg}$ is dense in $X$ (e.g. $Z=X_{\sing}$
or $Z=\bar S \setminus S$ in case $X=\bar S$ is the closure of a stratum in a Whitney stratification). 

We are interested in the intersection
\[ \bP(T^*_X M) \cap p^{-1}(\bP(\check{H}^{i})) = \bP(T^*_X M) \cap (M \times \bP(\check{H}^{i}))\]
for $\check{H}^{i}\subset \check{\bC}^N$ a generic linear subspace of codimension $i$ ($0\leq i \leq N-1$).
Note that $\bP(T^*_X M)$ is pure $(N-1)$-dimensional so that the dimension of this intersection
is bounded from below by $N-1-i$. For later use we also consider an additional proper algebraic subset
$Z'\subset \bP(T^*_X M)$ with $\dim \,Z'<N-1$, e.g. $Z'= \bP(T^*_X M)\cap \bP(T^*_{\bar{S'}} M)$ for $S'\neq S$
another stratum  in a Whitney stratification.\\

The natural action of the linear algebraic group $GL(\check{\bC}^N)$ on $\check{\bC}^N$ induces a transitive action on
$\bP(\check{\bC}^N)$ and on the  Grassmannian $\G^{i}(\check\bC^N)$ (or $\G^{i}(\bP(\check\bC^N))$)  
of (projective) linear subspaces of codimension $i$. 
So by  Kleiman's ``transversality result'' \cite{Kl}, one gets for {\em generic} 
$\check{H}^{i}$ the following result.

\begin{proposition}\label{p:polar}
For a given $i$, with $0\leq i \leq N-1$,
there is a Zariski-open dense set $\Omega^i(X)\subset \G^{i}(\check\bC^N)$ such that
the following properties are true for all $\check{H}^{i} \in \Omega^i(X)$:
\begin{enumerate}
\item $\bP(T^*_X M) \cap p^{-1}(\bP(\check{H}^{i}))$ is of pure dimension $N-1-i$, or empty.
\item $p^{-1}(\bP(\check{H}^{i}))$ intersects $\bP(T^*_U M)$ transversely.
\item $\bP(T^*_U M) \cap p^{-1}(\bP(\check{H}^{i}))$ is dense in $\bP(T^*_X M) \cap p^{-1}(\bP(\check{H}^{i}))$.
\item The dimension of $Z' \cap p^{-1}(\bP(\check{H}^{i}))$ is less than $N-1-i$.
\end{enumerate}
In particular, the intersection class
\[ [\bP(T^*_X M)] \cap [M \times \bP(\check{H}^{i})] = [\cl\bigl(\bP(T^*_U M) \cap p^{-1}(\bP(\check{H}^{i}))\bigr)]
\in CH_{N-1-i}(\bP(T^*M|_X))\]
equals $c^{1}(\cO(1))^{i}\cap [\bP(T_{X}^{*}M)]\in CH_{N-1-i}(\bP(T^*M|_X))$ and is therefore independent of
$\check{H}^{i} \in \Omega^i(X)$ (cf. \cite{Fu}).\fin
\end{proposition}

The Segre class, defined at (\ref{eq:segre}), is:
\[
s_{*}(T_{X}^{*}M) = \sum_{i\geq 0}\: \pi'_{*}(c^{1}(\cO(1))^{i}\cap [\bP(T_{X}^{*}M)]).
\]
In order to compute it, according to Proposition \ref{p:polar} we need to calculate the push down:
\[\hat{\pi}'_*([\bP(T^*_X M)] \cap [M \times \bP(\check{H}^{i})]) \in CH_{N-1-i}(X)\]
for $\hat{\pi}': \bP(T^*M|_X)\to X$ the (proper) projection. By definition and (a) above this is $0$ for
$k:=N-1-i>n$. Then let us assume $0\leq k \leq n$ and let us take a basis $\omega=(\omega_1, \ldots ,\omega_{k+1})$ of $\check{H}^{i}$,
where we identify the $\omega_j$ by the projection $p$ with the corresponding sections of $T^*M$ (i.e. the corresponding $1$-forms).
Then one has
\begin{equation*}
\begin{split}
\hat{\pi}'\bigl( \bP(T^*_U M) \cap p^{-1}(\bP(\check{H}^{i})) \bigr) &= \{x\in U|\: \rank\,(\omega_1, \ldots ,\omega_{k+1})\mid T_x U \leq k\}\\
&= \Crit (\omega_1, \ldots , \omega_{k+1})_U \:.
\end{split} \end{equation*}
So by (c) above we get
\begin{equation*}
\begin{split}
\hat{\pi}'\bigl( \bP(T^*_X M) \cap p^{-1}(\bP(\check{H}^{i})) \bigr) &= 
\overline{ \Crit (\omega_1, \ldots , \omega_{k+1})_U } \\
&= \overline{ \Crit (\omega_1, \ldots , \omega_{k+1})_{X_{\reg}} }
= P_{k}(X,\omega) \:
\end{split} \end{equation*}
where $P_{k}(X,\omega)$ is the \emph{polar locus of $\omega$}, see also Definition \ref{d:polar} for the special case $M=\bC^N$. Note that $\dim P_{k}(X,\omega)\leq k$ by (a) above.

For $k=0$ we have $ P_{0}(X,\omega)=\Crit (\omega_1)_{U}$, which is of dimension $\leq 0$ with
\[\hat{\pi}': \bP(T^*_U M) \cap p^{-1}(\bP(\check{H}^{i}))\to \Crit (\omega_1)_U\]
bijective so that
\[\hat{\pi}'_*([\bP(T^*_X M)] \cap [M \times \bP(\check{H}^{N-1})]) = [P_{0}(X,\omega)] \in CH_{0}(X) \:.\]
For the case $k>0$, 
we assume in addition now that the subspace $\check{H}^{i+1}\subset \check{H}^{i}$ with basis
$\omega':=(\omega_1, \ldots ,\omega_k)$ also belongs to the generic set $\Omega^{i+1}(X)$ of Proposition \ref{p:polar}.
Then $\dim P_{k-1}(X,\omega')< k$ and
\[\hat{\pi}': \bP(T^*_U M) \cap p^{-1}(\bP(\check{H}^{i}))\to \Crit (\omega_1, \ldots , \omega_{k+1})_U\]
is generically one to one so that $P_{k}(X,\omega)$ is also pure $k$-dimensional
or empty, with
\begin{equation} \label{eq:polarcycle}
\hat{\pi}'_*([\bP(T^*_X M)] \cap [M \times \bP(\check{H}^{i})]) = [P_{k}(X,\omega)] \in CH_{k}(X) \:.
\end{equation}

This finishes out proof, since in our context $T^*M$ is trivial and so $c^*(T^*M_{|X})=1$. 
\end{proof}

Let us set $P_{k}\bigl(X,(\omega_1, \ldots , \omega_{k+1})\bigr)$ $:= \emptyset$ for $k>n=\dim X$.
Now putting together
(\ref{eq:M=P}) and the commutative diagram (\ref{eq:conormal}),
 we get the following slightly more general version of Theorem \ref{t:MPcycle}.

\begin{theorem}\label{t:MPcycle2}
Let $X\subset M$ be a closed algebraic subset with $n:= \dim X<N$, endowed with a complex algebraic Whitney stratification
$\cS$, which has therefore finitely many strata $S$.
Assume $(\omega_1, \ldots , \omega_{N})$ is a general basis of $\check{\bC}^N$,
with $T^{*}M= M \times \check{\bC}^N$. Then
\[\sum_{S\in \cS}\; (-1)^{\dim S} \eta(S,\alpha) \cdot [P_{k}\bigl(\bar{S},(\omega_1, \ldots , \omega_{k+1})\bigr)]
= \check c_k^M(\alpha) \in CH_{k}(X)\]
for $0\leq k \leq n$ and all $\cS$-constructible functions $\alpha$.
\fin
\end{theorem}
The more precise meaning of ``general basis'' is that the subspace $\check{H}^{N-1-k}$ generated by $(\omega_1, \ldots , \omega_{k+1})$ 
belongs to  the generic set $\Omega^{N-1-k}(\bar{S})$ of Proposition \ref{p:polar}, for all
$0\leq k \leq n$ and all strata $S$, where $Z:=\bar S\setminus S$ and $U:=S$.

%%%
%%%%%%%%%%%%%%%%%%%%%%%%%%
\section{Degrees of affine polar varieties}\label{ss:deg}
%%%%%%%%%%%%%%
We now come back to the case when the ambient algebraic manifold $M$ is $\bC^N$,
with $X\subset \bC^N$ a proper algebraic subset of dimension $n<N$ which, for the moment, we assume
to be pure dimensional.
 Affine polar curves occured in the asymptotical equisingularity of families of hypersurfaces \cite{Ti-equi}. Degrees of affine polar varieties appeared implicitly in \cite{Ti-equi} and in the proof of the Lefschetz type formula for the global Euler obstruction \cite{STV},
which is the global counterpart of the L\^e-Teissier formula for the local Euler obstruction \cite{LT}.\\ 

To get the degrees $\gamma_k(X)$ of our polar varieties $P_k(X)$
and the generalized degrees $\gamma_k(\alpha)$ of our MacPherson cycles $\Lambda_k(\alpha)$,
we have to use additional genericity conditions. For this we look at the closure $\bar{X}$ of $X$ inside
the \emph{projective completion} $\bP^{N}$ of $M=\bC^N$ with $H^\ity:=\bP(\bC^N)\subset \bP^{N}$ 
the \emph{hyperplane at infinity}.
 We may endow $\bar X$
with an algebraic Whitney stratification $\hat\cS$ such that $\bar X \cap H^\ity$
is a union of strata and that its restriction to $X$ is the fixed Whitney stratification $\cS$ on $X$. 

By Kleiman's ``transversality result'' \cite{Kl}, one gets the following result, for {\em general} linear subspaces $\check{H}^{i}\subset \check{\bC}^N$ of codimension $i$.

\begin{proposition}\label{p:trans}
For a given $i$, with $0\leq i \leq N-1$,
there is a Zariski-open dense set $\hat\Omega^i(X)\subset \G^{i}(\check\bC^N)$ such that
the following property is true for all $\check{H}^{i} \in \hat\Omega^i(X)$:
\begin{equation}\label{eq:ker}
\bP(\ker \omega) \mbox{ is transversal inside } H^\ity \mbox{ to all strata 
of } \hat\cS \mbox{ which are included in } H^\ity. 
\end{equation}
Here $\bP(\ker \omega) \subset H^\ity$ is the projectivisation
of the $i$ dimensional linear subspace $\ker \omega\subset \bC^N$ for
$\omega:=(\omega_1, \ldots ,\omega_{k+1})$ a basis of $\check{H}^{i}$, where $k:=N-1-i$.
\fin
\end{proposition}
Note that $\ker \omega$ only depends on $\check{H}^{i}$, but not on the choice of the basis $\omega$.
This is not the case for the following result, where we use a basis 
$(\omega_1, \ldots ,\omega_N)$ of $\check{\bC}^N$ to simplify the notation.

\begin{proposition}\label{p:degree}
Let $\omega:=(\omega_1, \ldots ,\omega_{k+1})$, and for $k>0$ also $\omega':=(\omega_1, \ldots ,\omega_k)$ be chosen
so that  both $\bP(\ker \omega)$ and $\bP(\ker \omega')$ verify the condition (\ref{eq:ker}). Then the projection
 \[ (\omega_1, \ldots ,\omega_k) : P_k(X,\omega)\to \bC^k \]
is a proper and finite map. 
\end{proposition}

\begin{proof}
\smallskip
\noindent
Suppose that this map is not proper, where we first also assume $k>0$.
Then
there is a sequence of points $x_i\in X_\reg \cap P_k(X,\omega)$ which tends to
some point $y\in H^\ity$ with $\omega_j(x_i)$ bounded for all $1\leq j \leq k$.
Let $\hat\cS_\alpha\subset H^\ity$ be the stratum which contains $y$. 

Let us denote by $T_{x_i}$ the affine $n$-plane in $\bC^N$ tangent to $X_\reg$ at $x_i$ (which is obtained by translating  
the vector space $T_{x_i}X_\reg$ such that its origin becomes the point $x_i$). Let then 
$\langle T_{x_i} \rangle$ denote the projective closure of $T_{x_i}$ inside $\bP^{N}$,
i.e., $\langle T_{x_i} \rangle$ is the \emph{projective tangent space} of $X$ at $x_i$
(compare for example with \cite[p.181]{Har}).

 The assumed boundedness of $\omega_j(x_i)$ implies that $y\in \bP(\ker \omega')$.
Taking eventually a sub-sequence, we may assume without loss of generality that 
$\omega_j(x_i)$ converges to some point $b_j\in \bC$, for any $1\leq j \leq k$,
and that the limit 
of projective $n$-planes $\langle T_y \rangle:=\lim_i \langle T_{x_i}\rangle $ exists
inside the  Grassmannian $\G^{N-n}(\bP^N)$  
of projective linear subspaces of codimension $N-n$. 
Let $b:=(b_1, \ldots , b_k)$. Similarly, we denote by $\langle T_y\hat\cS_\alpha\rangle$
the projective plane in $H^\ity$ which concides with the tangent plane $T_y\hat\cS_\alpha$ 
at $y$, i.e., the projective tangent space of $\hat\cS_\alpha$ at $y$.

We claim that the assumed transversality $\bP(\ker \omega') \pitchfork_y \hat\cS_\alpha$ in $H^\ity$ 
implies that $\omega'$ is a submersion on  
 $T_{x_i}X_\reg$ for all $x_i$ close enough to $y$. If this were not true, then we 
would have the non-transversality in $\bP^N$, i.e.,
 \[ \overline{\omega'^{-1}(\omega'(x_i))} \not\pitchfork_{x_i}  \langle T_{x_i}\rangle .
\]
By passing to the limit, we get
\[ \overline{\omega'^{-1}(b)} \not\pitchfork_{y} \langle  T_y \rangle,
\]
which implies, since the plane $\overline{\omega'^{-1}(b)}$ is not included in $H^\ity$, that we have the non-transversality in $H^\ity$ of the sections
 \begin{equation}\label{eq:nontr_infty}
 (\overline{\omega'^{-1}(b)} \cap H^\ity) \not\pitchfork_{y}  (\langle T_y \rangle\cap H^\ity).
 \end{equation}

But we have $\overline{\omega'^{-1}(b)}\cap H^\ity = \bP(\ker \omega')$ and
$\langle T_y \rangle\supset \langle T_y\hat\cS_\alpha\rangle$ due to the Whitney (a)-regularity.
The last claim can be checked e.g. by pulling back under the canonical smooth projection $\pi: \bC^{N+1}\backslash \{0\}\to \bP^N$, with $\pi^{-1}(\bar{X})$ the affine cone over $\bar{X}$,
using the corresponding linear (tangent) spaces in $\bC^{N+1}$
(cf. \cite[p.183]{Har}) and the identification
$\G^{N-n}(\bP^N)\simeq \G^{N-n}(\bC^{N+1})$ of Grassmann manifolds.

Therefore (\ref{eq:nontr_infty})  contradicts the transversality $\hat\cS_\alpha \pitchfork \bP(\ker \omega')$.\\

By hypothesis, the linear map $\omega$ is not surjective on the tangent spaces $T_{x_i}X_\reg$. Due to the transversality of the intersection
$\overline{\omega'^{-1}(\omega'(x_i))} \cap \langle T_{x_i}\rangle $ proved above, this non-surjectivity is equivalent to the non-transversality:
\[  % \begin{equation}\label{eq:nontrans}
\overline{\omega_{k+1}^{-1}(\omega_{k+1}(x_i))} \not\pitchfork_{x_i} \langle T_{x_i} \rangle\cap \overline{\omega'^{-1}(\omega'(x_i))}\]
 in $\bP^N$.
Here we may also include the case $k=0$, with $\omega'$ defined as the zero map so that
$\overline{\omega'^{-1}(\omega'(x_i))}:=\bP^N$.
This non-transversality is then equivalent to the inclusion
\[ %\begin{equation}\label{eq:incl}
\overline{\omega_{k+1}^{-1}(\omega_{k+1}(x_i))} \supset \langle T_{x_i} \rangle\cap \overline{\omega'^{-1}(\omega'(x_i))},\]
since the left hand side is a hyperplane in $\bP^N$. Slicing both sides by $H^\ity$, we get
$ \bP (\ker \omega_{k+1}) \supset (\langle T_{x_i} \rangle\cap H^\ity) \cap \bP (\ker \omega')$.

Finally we pass to the limit, observing that only the space $\langle T_{x_i}\rangle$ varies, and we get
\[ \bP (\ker \omega_{k+1}) \supset ( \langle T_{y} \rangle\cap H^\ity) \cap \bP (\ker \omega').\]

This implies $y\in \bP(\ker \omega_{k+1})$ so that
$y\in \bP(\ker \omega)$ and contradicts the transversality 
$\hat\cS_\alpha \pitchfork \bP(\ker \omega)$ in $H^\ity$, since $\langle T_y \rangle\supset \langle T_y\hat\cS_\alpha\rangle$ by Whitney (a)-regularity.\\

So far we have shown the properness of the map $(\omega_1, \ldots ,\omega_k) : P_k(X,\omega)\to \bC^k$. Of course this map has then also finite fibres,
because $P_k(X,\omega)$ is an affine algebraic variety (or empty).
In particular, $\dim(P_k(X,\omega))\leq k$.
\end{proof}
%
%For $k=0$ we have that $P_0(X,\omega)$ is 
\begin{definition} \label{def:degree}
Assume that the linear subspace generated by  $\omega:=(\omega_1, \ldots ,\omega_{k+1})$,
and for $k>0$ also the linear subspace generated by  $\omega':=(\omega_1, \ldots ,\omega_k)$
belong to the generic sets of Propositions \ref{p:trans} and \ref{p:polar}
for a given $k$ with $0\leq k \leq n$. We have proved, for $k>0$, that the projection
 \[ (\omega_1, \ldots ,\omega_k) : P_k(X,\omega)\to \bC^k \]
is a proper and finite map, where $P_k(X,\omega)$ is pure $k$-dimensional or empty.
Its degree, which can be defined as the number of points in a general fibre,
will be denoted by $\gamma_{k}(X)$. 
In case $k=0$, we set $\gamma_{0}(X):= \# P_0(X,\omega)$, which makes sense since 
$P_0(X,\omega)$ is a finite set of points.
\end{definition}
The number $\gamma_{k}(X)$ will be independent of the choices of
the generic $\omega$, at least after restriction to another Zariski-open dense set of parameters.
This follows for example from the connectedness of the
Zariski-open dense sets, by using
(\ref{eq:mainindex}) and the fact that, for an algebraically constructible function $\alpha$
(like $\alpha= \Eu_X$), the correspondence
\[(\omega_1, \ldots, \omega_k, t_1, \ldots, t_k) \mapsto   \chi(X\cap \{\omega_1=t_1, \ldots, \omega_k=t_k\},\alpha)\]
defines an algebraically constructible function on the set of parameters (see e.g. \cite[sec.2.3]{Sch-book}). 
All the involved Zariski-open subsets depend
on the chosen Whitney stratification $\cS$ of $X$. Nevertheless the generic polar cycle classes
\[\check{c}^{\Ma}_{k}(X)= [P_k(X,\omega)] \in CH_{k}(X)\]
do not, and neither do the generic degrees $\gamma_k(X)$. One may notice that
 $\gamma_k(X)$ can also be described by
\[\gamma_k(X)= (\omega_1, \ldots ,\omega_k)_*([P_k(X,\omega)]) \in CH_{k}(\bC^k)\simeq \bZ .\]
%%%%%%%
%%%%%%%%%%%%%%%%%%%%%%%%%%%%%%%%%%%%%%%%%%%%%%%%%
\begin{note}\label{n:reduction-k=0}
The definition of $\gamma_k(X)$ in case $k>0$ can be reduced  to the case $k=0$ by taking iterated slices, as follows. Consider a general system of coordinates $\{x_j\}$ on $\bC^N$ (as discussed above) and  the inclusion
$i: X':=X \cap \{x_1=p_1, \dots, x_m=p_m\} \hookrightarrow X$ for general $p_1,\dots ,p_m$ (with $1\leq m\leq n$ fixed)
so that the affine space 
\[A_m:=\{x_1=p_1, \dots, x_m=p_m\}\subset \bC^N \] 
is transversal to the
given Whitney stratification of $X$. Consider the cartesian diagram for $m\leq k \leq n$
\begin{equation} \label{eq:intersec}
\begin{CD}
 \supp(\Lambda_k(\alpha)) @>>> X  @> (x_1,\dots,x_k) >>  \bC^k\\
@AAA  @A i AA @AAA \\
  \supp(\Lambda_k(\alpha))\cap A_m @>>> X'  @>  (x_{m-k+1},\dots,x_k) >> \bC^{k-m} @=  \{(p_1,\dots, p_m)\}\times \bC^{k-m}
\end{CD}
\end{equation}
From the refined Gysin homomorphism
$i^{!}$ (see \cite{Fu}) associated to this diagram, one gets  
the equalities
\begin{equation} \label{eq:gysin}
(-1)^m \cdot i^{!}\Lambda_k(\alpha) = \Lambda_{k-m}(\alpha|X') \quad \text{and} \quad
(-1)^m \cdot \gamma_k(\alpha) = \gamma_{k-m}(\alpha|X') \:.
\end{equation}
Here $i^!$ corresponds to the  intersection with $[A_m]$ inside $\bC^N$. The right hand side of the equalities make sense since the restriction $\alpha|X'$ is constructible with respect to the induced
Whitney stratification of $X'$ with strata given by the connected components of $S\cap A_m$.
Moreover the sign $(-1)^m$ is coming from 
\[ \dim(S\cap A_m)= \dim(S)-m\] 
for all strata $S$ and our sign convention in the definition 
(\ref{eq:macph}) of the $k$-th MacPherson cycle $\Lambda_k(\alpha)$. Of course we implicitly use
\begin{equation}\label{eq:restrNMD}
\eta(S,\alpha)= \eta(S\cap A_m,\alpha)
\end{equation}
for all strata $S$, since the corresponding ``normal Morse data'' do not change by taking general slices. In particular, $\eta(S\cap A_m,\alpha)$ depends only on $S$ and not on the choice of
a connected component of $S\cap A_m$.
For this it is enough to find a point $x$ in each connected component of $S\cap A_m$ such that $T^*_SM|_x$ contains a normally 
non-degenerate covector used for the definition of $\eta(S,\alpha)$ (as explained before), which then can also be used for the 
definition of $\eta(S\cap A_m,\alpha)$. The existence of such points follows for general $A_m$ from a simple dimension counting.
But if one only wants to assume that $A_m$ is transversal to our given stratification, then one has to use a result of Teissier
\cite[p.455]{Te-pol} telling us that 
the set of non-degenerate covectors is dense in all fibres of the projection
$T^*_SM\to S$, i.e., a complex analytic Whitney stratification is
$b_{\cod 1}$-regular, cf. \cite[Rem.5.1.9, p.320]{Sch-book}. 
Note that, in particular,
\[(-1)^m \cdot i^{!}\check c_*^M(\alpha)= \check c_*^M(\alpha|X') \:.\] 
This is the corresponding
{\em Verdier-Riemann-Roch Theorem} for the closed inclusion 
\[i: X' = X\cap A_m \to X \]
 of the transversal intersection
with the subspace $A_m\subset \bC^N$, whose normal bundle is trivial (cf. \cite[cor.0.1]{Sch4}).
\end{note}

\subsection{Interpretation of the degree and the index formula}\label{ss:index}
%%%%%%

We shall give a geometric interpretation of these degrees $\gamma_k(X)$, which will be essential in proving a Lefschetz type result
for affine pencils on $X$. In turn, this will provide a proof of Theorem \ref{t:mainindex}.
\begin{corollary}\label{c:degree}
Let $k$ be an integer, $0\leq k \leq n$.
Assume that the linear subspace generated by  $\omega:=(\omega_1, \ldots ,\omega_{k+1})$
and, for $k>0$, the linear subspace generated by  $\omega':=(\omega_1, \ldots ,\omega_k)$,
belong to the Zariski-open sets of Propositions  \ref{p:polar} and \ref{p:trans}.
Then the degree $\gamma_k(X)$ of the polar cycle $[P_k(X)]$ is the number of critical points of the pencil $\omega_{k+1}$ restricted to a general slice $X_\reg \cap \{ \omega_1 =t_1, \ldots ,\omega_k= t_k \}$, and all these critical points are complex stratified Morse points.
\end{corollary}
\begin{proof}
Let us first remark that for $k>0$ the generic choices of $\omega'$ and of the values
 $t_1, \ldots ,t_k$ imply that the plane $\{ \omega_1 =t_1, \ldots ,\omega_k= t_k \}$ intersects transversely the regular stratum $X_\reg$ and that the intersection points of
$P_k(X, \omega) \cap \{ \omega_1 =t_1, \ldots ,\omega_k= t_k \}$
are all contained in $X_\reg$ (compare with the discussion after Proposition \ref{p:polar}). 
This implies for $k=n+1$ that
\[ X\cap \{ \omega_1 =t_1, \ldots ,\omega_k= t_n\} = X_\reg \cap \{ \omega_1 =t_1, \ldots ,\omega_k= t_n \}\]
is  a finite set of points, which by definition are all complex Morse points for $\omega_{n+1}$.

For $k=0$ this is by definition just $P_0(X, \omega)\subset X_\reg$.
These intersection points are singular points of the pencil $\omega_{k+1}$ on the slice 
$X_\reg \cap \{ \omega_1 =t_1, \ldots ,\omega_k= t_k \}$. They are complex  Morse points 
by the transversality result (b) of Proposition \ref{p:polar}.
In fact, for $k=0$, this just means in local coordinates that the graph of the $1$-form $\omega_1$ intersects $T^*_{X_{\reg}}M$ 
transversely inside $T^*M$, which corresponds to a complex Morse point.
If we work with an ambient algebraic Whitney stratification, then we even can get, 
by the dimension estimate (d) of Proposition \ref{p:polar}, that these are stratified Morse points.
 The case $k>0$ is reduced to this by slicing with $\{ \omega_1 =t_1, \ldots ,\omega_k= t_k \}$.
\end{proof}

We continue to work with some fixed algebraic Whitney stratification $\cS$ of $X \subset \bC^N$ and some extension of it to a stratification $\hat \cS$ of $\bar X$.
 Then the polar classes $[P_{k}(\bar S)]$ of any stratum $S \in \cS$,
as well as their degrees $\gamma_k(\bar S)$, are well defined. In particular, we have $P_{k}(\bar S) = \bar S$ for $k = \dim S$ and by definition $P_{k}(\bar S)= \emptyset$ for $k> \dim S$ so that
\[  \gamma_{k}(\bar S)= 0 \quad \text{for $k > \dim S$ and} \quad \gamma_{\dim S}(\bar S)= \deg S. \]

 For $k=N-1-i$ with $0\leq k \leq n$, assume that the linear subspace generated by  $\omega:=(\omega_1, \ldots ,\omega_{k+1})$
and, for $k>0$,  the linear subspace generated by  $\omega':=(\omega_1, \ldots ,\omega_k)$,
belong to the generic sets $\hat \Omega^i(\bar S)$ and $\Omega^i(\bar S)$ of Propositions \ref{p:polar} and \ref{p:trans}, for all strata $S$.
This is a finite intersection of Zariski-open dense subsets, hence it is also a Zariski-open dense subset.

\begin{proposition}\label{p:lefschetz}
Let $\omega$ be chosen as explained before, with $(t_1, \dots, t_k)$ generic values of $\omega'$.
 Then, for general $t_{k+1}$, the slice  $X\cap \{\omega_1=t_1, \ldots, \omega_k=t_k\}$
is  homeomorphic to the space obtained from the slice  $X\cap \{\omega_1=t_1, \ldots, \omega_{k+1}=t_{k+1}\}$
by attaching to it the local Morse data of the stratified complex Morse points of the pencil $\omega_{k+1}$.
\end{proposition}
\begin{proof}
The generic choices of $\omega$, of its generators $\omega_1, \ldots , \omega_k$ and of the values
 $t_1, \ldots ,t_k$ imply that the plane $\{ \omega_1 =t_1, \ldots ,\omega_k= t_k \}$ 
 intersects transversely all the strata of the stratification $\hat \cS$.
 In particular, it does not intersect any stratum of dimension less than $k$.
 In this way, the slice $X' := X \cap \{\omega_1=t_1, \ldots, \omega_k=t_k\}$
is endowed with the induced Whitney stratification. 

It follows from the definition of $\hat \Omega^i$ that the axis $\bP(\ker \omega)$ of the affine pencil defined by $\omega_{k+1}$ is transversal to the strata of the slice $X'$. 
Therefore the pencil $\omega_{k+1}$ is a locally trivial fibration in the neighbourhood
of $H^\ity$ and the hyperplanes of the pencil intersect transversely
all the strata of the stratification $\cS$, except at a finite number of points
on each stratum. These points are the intersections of $X'$ with the polar loci
$P_k(\bar S)$, for $S\in \cS$. As shown in Corollary \ref{c:degree},
all such points are complex stratified Morse points of the pencil $\omega_{k+1}$.

 By the Lefschetz theory and stratified Morse theory \cite{GM}, the total space $X'$ of the pencil is obtained from the general
 hyperplane section $X'\cap \{ \omega_{k+1} = t_{k+1}\}$ by attaching 
the local Morse data of the stratified critical points
 of the pencil $\omega_{k+1}$, which are complex stratified Morse points in our case.

 Note that this argument also works for $k=n+1$ with $X'\cap \{ \omega_{n+1} = t_{n+1}\}=\emptyset$ and
 $X'$ a finite set.
\end{proof}

\subsection{Proof of Theorem \ref{t:mainindex}.} 
We apply the Euler characteristic $\chi(\cdot, \alpha)$ weighted by the constructible function $\alpha$ to the decomposition provided by Proposition \ref{p:lefschetz}.
 We start with generic choices of $\omega, \omega'$, of its generators $\omega_1, \ldots ,\omega_k$ and of the values $t_1, \ldots, t_k$
 as before.
The number of singular points of the pencil $\omega_{k+1}$
on every stratum 
\[S':=S\cap \{\omega_1=t_1, \ldots, \omega_k=t_k\}\]
 is equal to $\gamma_{k}(\bar S)$, by Corollary \ref{c:degree} for $\dim S\geq k$,
and $S'=\emptyset$ for  $\dim S< k$. 
Note that here we allow a stratum $S'$ to be disconnected.
We get:
\[ \chi(X\cap \{ \omega_1= t_1, \ldots ,\omega_{k}= t_{k} \}, \alpha) -
\chi(X\cap \{ \omega_1= t_1, \ldots ,\omega_{k+1}= t_{k+1} \}, \alpha)  \]
\[ =
\sum_{S\in \cS} 
\gamma_{k}(\bar S) \chi(\LMD(\omega_{k+1}|S', \alpha)),
\]
where $\LMD(\omega_{k+1}|S')$ denotes the local Morse datum of $\omega_{k+1}$ in its critical points on $S'$
considered as a stratum of $X\cap \{ \omega_1= t_1, \ldots ,\omega_{k}= t_{k} \}$.
They are all isomorphic for different critical points in $S'$, because they are
all complex stratified Morse critical points,
and the corresponding normal Morse data are all isomorphic to the normal Morse datum of the stratum
$S$. Notice that the terms of the sum corresponding to strata $S$ of dimension less than $k$
are zero since $\gamma_{k}(\bar S)=0$ in that case.

 Finally we use the decomposition 
 of the local Morse datum $\LMD(\omega_{k+1}|S')$, 
 which was proven in the stratified Morse theory by Goresky and MacPherson \cite{GM}, as a product of normal and tangential Morse data.
 At the level of the Euler characteristic, we get, using Lemma \ref{thm:indmain},
\[ \chi(\LMD(\omega_{k+1}|S', \alpha)) = (-1)^{\dim S'} \cdot
\chi (\NMD(S',\alpha)) \:,\]
with $\dim S' = \dim S-k$ and
\[\chi (\NMD(S',\alpha)) = \chi (\NMD(S),\alpha)) =: \eta(S,\alpha)\]
 the normal index.
Note that the first equality is explained in (\ref{eq:restrNMD}).  
Now formula (\ref{eq:mainindex}) follows from the definition 
\[\gamma_k(\alpha) := \sum_{S\in \cS} (-1)^{\dim S}\cdot  \eta(S,\alpha)\gamma_k(\bar S) \:.\]

%%%%%%%%%%%%%%%%%%%%%%%%%%%%%%%%%%%%%%%%%%%%%%%%%%%%%%%

\end{document}